%-----------------------------------------------------------------------------
% 
%Beginning of Quadratic minima and modular forms II.tex
%-----------------------------------------------------------------------------
%
% 
%AMS-TeX 2.1 file for journals.
%
% 
%
\input amstex
\documentstyle{amsppt}
\NoBlackBoxes

\topmatter
\title Quadratic minima and modular forms II\endtitle
\author Barry Brent\endauthor

\address 1302 Eleanor Avenue, St. Paul, MN 55116 \endaddress

\email barryb\@primenet.com\endemail


\abstract We give upper bounds
on the size of the gap between a non-zero constant term and the next
 non-zero Fourier 
coefficient
of an entire modular form
for $\Gamma_0(2)$.
We derive upper bounds 
for the minimum positive integer represented by
level two even positive-definite quadratic forms. 
The bounds are sharper than the ones
we proved, and 
slightly sharper than the ones we conjectured, in part I. 
\endabstract

\thanks The author is grateful to his advisor, Glenn Stevens, for 
numerous helpful conversations.
\endthanks


\endtopmatter

\document

\head 1. Introduction\endhead

Carl Ludwig Siegel showed in 
[Siegel 1969]  (English translation, [Siegel 1980])
that the constant terms of certain level one negative-weight 
modular forms $T_h$ are non-vanishing (`` \it Satz 2 \rm ''), 
and that this
implies an upper bound on the least positive exponent of a 
non-zero Fourier 
coefficient for  any level one entire modular form of weight 
$h$ with a 
non-zero constant term. 
Level one theta functions fall into this category. 
Their Fourier coefficients code up
representation numbers 
of quadratic forms. 
For positive even $h$, Siegel's result gives an 
upper bound on the
least positive integer represented by a
 positive-definite even unimodular quadratic form in $n = 2h$ variables. 
This bound is sharper than Minkowski's for large $n$. (Mallows, Odlyzko and 
Sloane have improved Siegel's bound in [Mallows, Odlyzko, and Sloane 1975].)

 John Hsia [private communication
to Glenn Stevens] 
suggested that Siegel's approach might be extended to higher levels. 
Following this hint, we 
constructed an analogue  of $T_h$ for $\Gamma_0(2)$,
which we denote as $T_{2,h}$. To prove $Satz$ 2, Siegel controlled the 
sign of the 
Fourier coefficients in the principal part of $T_h$. 
In [Brent 1998] (henceforth, ``part I''),
following Siegel, 
we 
found upper bounds
 for the first positive exponent of a non-zero
Fourier coefficient occuring in the expansion at infinity 
of an entire modular form with a non-zero constant term for $\Gamma_0(2)$
in the case \newline $h \equiv 0$ $(\hskip -.07in \mod 4)$. 
Siegel's method carried over intact. 

In part I, we also stated that
it was not clear that Siegel's method forces the non-vanishing of 
the $T_{2,h}$ 
constant terms when $h \equiv 2$ $(\hskip -.07in \mod 4)$. 
But it turns out that we can tweak our definition 
of the $T_{2,h}$ and carry out Siegel's strategy.

Let us denote the vector space of entire modular forms of weight $h$ 
for $\Gamma_0(N)$ as $M(N,h)$. In part I,
we proved that the second non-zero Fourier coefficient
of an 
an element of  $M(2,h)$ with non-zero constant term must have exponent
at most $\roman{dim} M(2,h) = 1+\left \lfloor \frac h4 \right \rfloor
= r$ (say) if $h \equiv 0 \,(\hskip -.1in \mod 4)$.
This corresponds exactly to Siegel's bound for $f \in M(1,h).$
If  $h \equiv 2 \,(\hskip -.1in \mod 4)$, 
however, we only showed that 
the exponent is no more than $2r.$ 

In the second section,  we prove that the exponent is
at most $r$ if $h \equiv 2 \,(\hskip -.1in \mod 4)$. 
In the third section, we apply this result to the theory of quadratic forms.
We show that, if $Q$ is an even positive-definite level two quadratic form 
in $v=8u+4$ variables, then $Q$ represents
a positive integer $2n \leq 1+\frac{v}4$. (In part I, we  obtained the
weaker bound $ 2 + \frac v2.$ We also showed that
if $v=8u,$ then $Q$ represents an even positive integer $\leq 2+\frac v4.$)

\head 2. Bounds for gaps in the Fourier expansions of 
entire modular forms \endhead
Section 2.1 is introductory. All but one of
the results are stated without proof.
The reader is referred to part I for
details.
In section 2.2, we estimate the first positive exponent of a non-zero Fourier
coefficient in the expansion of an entire modular form for 
$\Gamma_0(2)$ with a non-zero constant term.

\subhead {\rm 2.1.}\enspace Some modular objects\endsubhead 
This section is a tour of the objects mentioned in the article. 
The main building blocks are 
Eisenstein series with known divisors and computable Fourier expansions. 

As usual, we denote by 
$\Gamma_0 (N)$ the congruence subgroup 

$$
\Gamma_0 (N) = \left\{\left(\matrix a & b\\c & d\endmatrix \right)\in\, 
SL(2,\bold Z): 
c \equiv 0\,(\mod N )\right\}.
$$

The vector space of entire modular forms of one variable in the upper half plane 
$\goth{H}$ of weight $h$  
for $\Gamma_0 (N)$ (``level $N$'') and trivial character, we denote by $M(N,h)$. 
We have an inclusion lattice satisfying:
$$ M(L,h) \subset M(N,h) \roman{ \,\,\,\,if\, and\, only\, if\, } L | N.$$
More particularly, any entire modular form for $SL(2,\bold Z)$ is also one for 
$\Gamma_0(2).$

The dimension 
of $M(N,h)$ is denoted by $r(N,h)$, or $r_{_h}$, or by $r$. 
For any positive even $h$, 
$$ r(2,h) = \left \lfloor\frac{h}{4} \right \rfloor+1.$$

We write $\Delta$ for the weight 12, level one cusp form 
with Fourier series
$$
\Delta = \sum_{n=1}^{\infty}\tau(n)\,q^{\,n} \,.
$$
%%%%%%%%%%
and product expansion
$$
\Delta = q\prod_{n=1}^{\infty} (1-q^n)^{24} \,.
$$
Here, $\tau$ is the Ramanujan function. 

We describe some level two objects, using three special divisor sums:
$$\sigma^{\roman{odd}}(n) = \sum \Sb 0<d|n\\d\,\,\roman{odd}\endSb d,$$ 
$$\sigma_k^{\roman{alt}}(n) = 
\sum_{0<d|n } (-1)^{d}d^k,$$\flushpar and $$\sigma_{N,k}^{*}(n) =
\sum_{\Sb  0<d|n\\N \not\ \hskip -.01in | \frac nd 
\endSb } d^k.$$ \flushpar
Let $E_{\gamma,2}$ denote the unique normalized form in 
the one-dimensional space $M (2 , 2)$
 (\it i.e.\rm \,the leading coefficient in the Fourier expansion of the form is a $1$). 
The Fourier series is
$$
E_{\gamma,2} = 1+24\sum_{n=1}^{\infty} \sigma^{\roman{odd}}(n)\,q^{\,n}.\tag{2\,-1}
$$
\flushpar
%%%%%%%%%%%%%%%%%%%%
$E_{\gamma,2}$ has a $\frac12$-order zero at points of 
$\goth{H}$ which are $\Gamma_0 (2)$\,-equivalent to  
$-\frac12 +  \frac12 i = \gamma$ (say).
The vector space $M (2,4)$ is spanned by two forms $E_{0,4}$ and $E_{\infty,4}$, 
which vanish with order one at the $\Gamma_0 (2)$\,-inequivalent zero and infinity cusps, 
respectively. They have Fourier expansions
$$
 			E_{0,4}  = 1 + 16\sum_{n=1}^{\infty}\sigma_3^{\roman{alt}}(n)\,q^{\,n}\tag{2-2}
$$
\flushpar
and
$$
			E_{\infty,4}  = \sum_{n=1}^{\infty}\sigma_{2,3}^{*}(n) \,q^{\,n}.\tag{2-3}
$$

The following lemma was proposition 2.3 in part I:

\proclaim{Lemma 2.1} The modular form $E_{\infty,4} \in M(2,4)$ has the following
product decomposition in the variable $q$:
$$
E_{\infty,4}(z) = 
q\prod_{0<n\in 2\bold Z} \left(1-q^n \right)^8 \prod_{0<n\in \bold Z\backslash 2\bold Z} 
\left(1-q^n \right)^{-8}.
$$
\endproclaim

We don't need it, but it is also easy to show that

$$
E_{0,4}(z) = 
q\prod_{0<n\in 2\bold Z} \left(1-q^n \right)^8 \prod_{0<n\in \bold Z\backslash 2\bold Z} 
\left(1-q^n \right)^{16}.
$$
To exploit Lemma 2.1, we need a result due essentially
to Euler (lemma 2.11 of part I, quoted from
 [Apostol 1976], Theorem 14.8, valid when the infinite series is 
absolutely convergent):
\vskip .2in
\proclaim{Lemma 2.2} For a given set A and a given arithmetical function
f, the numbers $p_{A,f}(n)$ defined  
by the equation 
$$\prod_{n \in A}(1-x^n)^{-f(n)/n}\, = 1 + \sum_{n=1}^{\infty}p_{A,f}(n)x^n$$
satisfy the recursion formula $$np_{A,f}(n)\,= \sum_{k=1}^n f_A(k)p_{A,f}(n-k),   $$
where $p_{A,f}(0)=1$ and $$f_A(k) = \sum_{\Sb  d|k\\d \in A \endSb }f(d).$$

\endproclaim

Next, we construct a level two analogue of the level one Klein invariant $j$:
%%%%%%%%%%%%%%%%%%%%%%%%%%%%%%%%%%%%%%%%%%%%%%%%%%%%%%%%%%%%%%%%%%%%%%%%%%%%%%%%
$$
			j_{_2} =\, E_{\gamma,2}^2E_{\infty,4}^{-1} .
$$
The function $j_{_2}$ is analogous to $j$ because it is
 modular (weight zero) 
for 
$\Gamma_0 (2)$, holomorphic on 
the upper half plane, has  a
simple pole at infinity, generates the 
field of $\Gamma_0 (2)$\,-modular functions, and 
defines a 
bijection of a  $\Gamma_0 (2)$ fundamental set with \bf C\rm.

The following lemma was proposition 2.7 in part I:

\proclaim{Lemma 2.3}
 For $z \in \goth{H}$, 
$$ \frac d{dz}\,j_2(z) = -2 \pi i E_{\gamma,2}(z) E_{0,4}(z) E_{\infty,4} (z)^{-1}. $$
\endproclaim

We introduce a 
level two analogue of Siegel's $T_h$ for 
$h\equiv 2 \,(\hskip -.1in \mod 4)$.  For $r  = r (2,h)$,
we set

$$
			t_{h} =E_{0,4}E_{\infty,4}^{-r}.
$$
 
\flushpar The $t_h$ are a useful
replacement for the functions
$$
			T_{2,h} =E_{\gamma,2}^2E_{0,4}E_{\infty,4}^{-1-r}
$$
\flushpar we defined in part I for the same $h$.

  Finally, for $h \equiv 2 \,(\hskip -.07in\mod 4)$ and
$f \in M(2,h)$, let $$w_2(f) = E_{\gamma,2}^{-1}  E_{\infty,4}^{1-r}f.$$ 
This replaces the part I function 
 $W_2(f) = E_{\gamma,2} E_{\infty,4}^{-(h+2)/4}f$ defined on the same $f$.

\subhead {\rm 2.2.}\enspace  A structural result on the $q$ series of 
entire level two modular forms\endsubhead 

We establish a sequence of propositions  mimicking the argument of
[Siegel 1980], pp. 249-254. (Siegel's proof is also sketched
on p. 263 of part I.)

\proclaim{Proposition 2.1} The map $w_2$ 
 is a vector space isomorphism from $M(2,h)$
onto the space of polynomials in $j_2$ of degree less than $r$.
\endproclaim

\flushpar \it Proof. \rm  No
 non-trivial polynomial in $j_2$ 
can vanish almost everywhere, so the modular forms
$j_2^d E_{\gamma,2} E_{\infty,4}^{r-1}$, $d = 0, 1, ..., r-1$ 
 are a basis for $M(2,h),$    and we have
$$w_2(j_2^d E_{\gamma,2} E_{\infty,4}^{r-1}) = j_2^{d}.$$ The map $w_2$ 
is clearly linear and 1-to-1. $\hskip .1in \boxed{} $

\proclaim{Proposition 2.2} For $f \in M(2,h)$,
the constant term in the Fourier expansion at infinity of $t_{h}f$ is
zero. 
\endproclaim

\flushpar \it Proof. \rm By applying Lemma 2.3, we see that
 $ w_2(f) \frac d{dz} j_2 = -2 \pi i  t_{h}f.$  Thus, 
$t_{h}f$ is the derivative of a polynomial in $j_2$, so it
can be expressed in a neighborhood of infinity as the derivative 
with respect to $z$ of a power series in the variable $q = \exp (2 \pi i z)$.
This derivative is a power series in $q$ with vanishing constant term.
 $\hskip .1in \boxed{} $

%%%%%%%%%%%%%%%%%%%%%%%%%%%%%%%%%%%%%%%%%%%%%%%%%%%%%%%%%%%%%%%%%%%%%%%%%%%%%%%%

\proclaim{Proposition 2.3} For $h \equiv 2 \, (\hskip -.07in 
\mod 4)$, the constant term in the Fourier expansion
at infinity of $t_{h}$ is non-zero.
\endproclaim
 
\flushpar \it Proof. \rm Lemmas 2.1 and 2.3
imply that, for fixed $s$,

$$
E_{\infty,4}^{-s}=q^{-s}\sum_{n=0}^{\infty} R(n)q^n,\tag{2-4}
$$ 

\flushpar
where $R(0)=1$ and $n>0$ implies that
$$
R(n)=\frac{8s}n\sum_{a=1}^n \sigma_1^{\roman{alt}}(a)R(n-a). \tag{2-5}
$$
 The divisor functions
$\sigma_k^{\roman{alt}}(n)$, $k$ odd, alternate sign, 
so the alternation of the sign of $R(n)$ 
follows by an easy induction argument from (2-5). To be specific, 
$R(n)= U_n(-1)^n$ for some $U_n > 0$. Thus we may write
$$ E_{\infty,4}^{-r} = U_0(-1)^0 q^{-r} + U_1 (-1)^1 q^{1-r} + ... +
U_{r-1}(-1)^{r-1}q^{-1} + U_r(-1)^r + ... \, .
$$  \newline

On the other hand, the Fourier coefficient of $q^n,\,\, n \geq 0$, in the expansion of 
$E_{0,4}$
is $W_n(-1)^n$ for positive $W_n$, by (2-2). Thus the constant term of 
$ t_h = E_{0,4} 
E_{\infty,4}^{-r}$
is $$\sum_{m=-r}^0 U_m (-1)^{m}W_{r-m} (-1)^{r-m} \neq 0.
$$ 
$\hskip 4in \boxed{} $
\vskip .1in

What follows is our main theorem on modular forms.

\proclaim{Theorem 2.1}
Suppose $f \in M(2,h)$ with Fourier expansion at infinity
$$f(z) = \sum_{n=0}^{\infty} A_n q^n, \,\,A_0 \neq 0.$$
Then some $A_n \neq 0, 1 \leq n \leq r(2,h)$.

\endproclaim

\flushpar \it Proof. \rm First suppose that 
 $h \equiv 2 \, (\hskip -.07in \mod 4)$.
We  denote the coefficient of $q^n$ in the Fourier expansion of $f$ at infinity
as $c_n[f].$ 
The normalized meromorphic form $t_{h}$  has a Fourier 
series of the form
$$t_h  = C_{h,-r}q^{-r} + ... + C_{h,0} + ... ,$$
with $C_{h,-r}  = 1$. By Proposition 2.2, 
$$0 = c_{_0}[t_{h} f ] = C_{h,0}A_0 + ... + C_{h,-r} A_r.$$ 
By hypothesis, $A_0 \neq 0$. By Proposition 2.3, $C_{h ,0} \neq 0$, so 
$$	A_0 = -(C_{h ,0})^{-1}(C_{h ,-1}A_1 + ... + C_{h ,-r} A_r  ).$$
It follows that one of the $A_n  \,(n = 1, ..., r ) $ is non-zero.

To complete the proof, we point out that the claim was
proved for $h \equiv 0 \, (\hskip -.07in \mod 4)$
in Theorem 2.12 of part I by the same sort of argument.

 $\hskip 4in \boxed{} $
\vskip .1in
\it Remark\rm \,\,This result is slightly better than the bound $n \leq r+1$ 
of conjecture 6.2, part I for  $h \equiv 2 \, (\hskip -.07in \mod 4)$. 
The reason is our different choice
of a level two $T_h$ analogue.
%%%%%%%%%%%%%%%%%%%%%%%%%%%%%%%%%%%%%%%%%%%%%%%%%%%%%%%%%%%%%%%%%%%%%%%%%%%%%%%
\head 3. Quadratic minima \endhead
\subhead{\rm 3.1} Quadratic forms and modular forms \endsubhead
For even $v$, set $\bold{x} = \, ^t(x_1,\,...,\, x_v)$, so that
$\bold{x}$ is a column vector. Let $A$ be an $v$ by $v$ square symmetric matrix 
with integer entries, even entries on the diagonal, and positive eigenvalues.
 Then $Q_A(\bold{x}) = \, ^t\bold{x} A \bold{x}$ is a homogenous second degree
polynomial in the $x_i$. We refer to $Q_A$ as the 
even positive-definite quadratic form associated to $A$.
If $\bold{x} \in \bold{Z}^v$, then $Q_A(\bold{x})$ is a non-negative even number,
which is zero only if $\bold{x}$ is the zero vector. The level of $Q_A$
is the smallest positive integer $N$ such that $NA^{-1}$ also has integer entries
and even entries on the diagonal.
 Let $\#Q_A^{-1}(n)$ denote 
the cardinality of
the inverse image in $\bold{Z}^v$ of an integer $n$ under the quadratic form $Q_A$.

The following specialization of known
results was proposition 5.1 of part I.

\proclaim{Lemma 3.1} Suppose that $Q_A$ is a level two quadratic form.
Then the function 
$\Theta_A: \goth{H} \rightarrow \bold{C}$ satisfying
$$\Theta_A(z) = \sum_{n=0}^{\infty} \#Q_A^{-1}(2n)q^n\tag{5-1}$$
lies in $M(2,\frac v2)$.\endproclaim

Since $M(2,h)$ is non-trivial only for even $h$, it also follows that $4 | v$.
\subhead{\rm 3.2} Quadratic minima \endsubhead
In this section we apply Theorem 2.1 to the problem of quadratic minima.

\proclaim{Theorem 3.1} If $Q$ is a level two  
even positive-definite quadratic form
in $v$ variables, $8 | v$, then $Q$ represents
a positive integer $2n \leq 2 + \frac v4$. If $v \equiv 4$
 $(\hskip -.07in \mod 8)$, then $Q$ represents a positive integer 
$2n \leq 1 + \frac v4$. 
\endproclaim

\flushpar \it Proof. \rm The claim for $8|v$ was included in
theorem 5.2 of part I. Suppose
 $v=8u+4$. 
Let $A$ be the matrix associated to $Q$, so that $Q=Q_A$. 
Then $\Theta_A \in M(2,4u+2)$,
and $\#Q_A^{-1}(2n) \neq 0$ for some $n, \, 1 \leq n \leq r(2,4u+2) = 1 + u$. Thus
$Q$ represents
an integer $2n \leq 2 + 2u = 1 + \frac v4$.$\hskip .1in \boxed{} $

\refstyle{C}

\enddocument